\documentclass[reqno,a4paper]{amsart}%[rm]{birkmult}%
\usepackage{amssymb}
\usepackage[colorlinks=true,hypertex]{hyperref}%,pagebackref
\allowdisplaybreaks[4]
\theoremstyle{plain}
\newtheorem{thm}{Theorem}

\theoremstyle{remark}
\newtheorem{rem}{Remark}
\DeclareMathOperator{\td}{d\mspace{-2mu}}
\date{Commenced on 29 January 2009 and completed on 31 January 2009 in Melbourne}
\date{}

\begin{document}

\title[Oppenheim's inequality relating to the cosine and sine functions]
{A concise proof of Oppenheim's double inequality relating to the cosine and sine functions}

\author[F. Qi]{Feng Qi}
\address[F. Qi]{Research Institute of Mathematical Inequality Theory, Henan Polytechnic University, Jiaozuo City, Henan Province, 454010, China}
\email{\href{mailto: F. Qi <qifeng618@gmail.com>}{qifeng618@gmail.com}, \href{mailto: F. Qi <qifeng618@hotmail.com>}{qifeng618@hotmail.com}, \href{mailto: F. Qi <qifeng618@qq.com>}{qifeng618@qq.com}}
\urladdr{\url{http://qifeng618.spaces.live.com}}

\author[B.-N. Guo]{Bai-Ni Guo}
\address[B.-N. Guo]{School of Mathematics and Informatics, Henan Polytechnic University, Jiaozuo City, Henan Province, 454010, China}
\email{\href{mailto: B.-N. Guo <bai.ni.guo@gmail.com>}{bai.ni.guo@gmail.com}, \href{mailto: B.-N. Guo <bai.ni.guo@hotmail.com>}{bai.ni.guo@hotmail.com}}
\urladdr{\url{http://guobaini.spaces.live.com}}

\begin{abstract}
In this paper, we provide a concise proof of Oppenheim's double inequality relating to the cosine and sine functions. In passing, we survey this topic.
\end{abstract}



\thanks{The first author was partially supported by the China Scholarship Council}

\thanks{This paper was typeset using \AmS-\LaTeX}

\maketitle

\section{Introduction and main results}

In \cite{Oppenheim-monthly-57}, the following problem was posed: For each $p>0$ there is a greatest $q$ and a least $r$ such that
\begin{equation}\label{Oppenheim-monthly-57-ineq}
\frac{q\sin x}{1+p\cos x}\le x\le\frac{r\sin x}{1+p\cos x}
\end{equation}
for $0\le x\le\frac\pi2$. Determine $q$ and $r$ as functions of $p$.
\par
In \cite{Carver-Monthly-58}, it was explicitly obtained that
\begin{enumerate}
\item
the least value of $r$ required by the problem is
\begin{equation*}
\begin{aligned}
r&=\frac\pi2 & & \text{when $p\le\frac\pi2-1$,}\\
r&=p+1 & & \text{when $p\ge\frac\pi2-1$;}
\end{aligned}
\end{equation*}
\item
the required greatest value of $q$ is
\begin{equation*}
\begin{aligned}
q&=p+1 & & \text{when $p\le\frac12$,}\\
q&=\frac\pi2 & & \text{when $p\ge\frac\pi2$.}
\end{aligned}
\end{equation*}
\end{enumerate}
\par
In~\cite[p.~238, 3.4.15]{mit}, it was listed that
\begin{equation}
\frac{(p+1)\sin x}{1+p\cos x}\le x\le \frac{(\pi/2)\sin x}{1+p\cos x}
\end{equation}
for $0<p\le\frac12$ and $0\le x\le\frac\pi2$.
\par
In~\cite[p.~521, (26)]{construct}, by Techebysheff's integral inequality, it was constructed that
\begin{equation}\label{cusa-2}
\frac{\sin x}{x}\ge\frac{1+\cos x}{2},\quad 0<x\le\frac\pi2
\end{equation}
and
\begin{equation}
\frac{\sin x}{x}\ge\frac{1+2\cos x}{3}+\frac{x\sin x}6,\quad 0<x\le\frac\pi2.
\end{equation}
The inequality~\eqref{cusa-2} can be rewritten as
\begin{equation}\label{cusa-2-rewrit}
\frac{2\sin x}{1+\cos x}\ge x,\quad 0\le x\le\frac\pi2.
\end{equation}
\par
In \cite{sandor-bencze}, it was pointed out that the inequality
\begin{equation}\label{Cusa's-ineq}
\frac{3\sin x}{2+\cos x}<x
\end{equation}
was discovered by Nicolaus de Cusa (1401--1464) using certain geometrical constructions. In~\cite{sandor-Erdelyi-02}, the inequality~\eqref{Cusa's-ineq} was generalized as follows: For $a,b,c>0$ such that $2b\le c\le a+b$,
\begin{equation}
\frac{c\sin x}{a+b\cos x}<x,\quad 0<x<\frac\pi2.
\end{equation}
This is equivalent to the left-hand side inequality in~\eqref{Oppenheim-monthly-57-ineq} for $2p\le q\le1+p$.
\par
In \cite[Theorem~7]{Oppeheim-Zhu-mia}, a complete answer to the above problem was obtained as follows: Let $0\le x\le\frac\pi2$ and $p>0$, then the inequality~\eqref{Oppenheim-monthly-57-ineq} holds in cases:
\begin{enumerate}
\item
When $0<p<\frac12$, we have $q=p+1$, $r=\frac\pi2$;
\item
When $\frac12\le p<\frac\pi2-1$, we have $q=4p(1-p^2)$, $r=\frac\pi2$;
\item
When $\frac\pi2-1\le p<\frac2\pi$, we have $q=4p(1-p^2)$, $r=p+1$;
\item
When $\frac2\pi\le p<\infty$, we have $q=\frac\pi2$, $r=p+1$.
\end{enumerate}
\par
The aim of this paper is to provide a concise proof of the inequality~\eqref{Oppenheim-monthly-57-ineq}.
\par
Our main results may be recited as the following theorems.

\begin{thm}\label{Opeheim-Sin-Cos.tex-thm-1}
For $p>0$ and $x\in\bigl(0,\frac\pi2\bigr]$, let
\begin{equation}
f_p(x)=\frac{\sin x}{x(1+p\cos x)}.
\end{equation}
Then the function $f_p(x)$ is strictly
\begin{enumerate}
\item
increasing when $p\ge\frac2\pi$;
\item
decreasing when $p\le\frac12$.
\end{enumerate}
Moreover, when $\frac12<p<\frac2\pi$, the function $f_p(x)$ has a unique maximum on $\bigl(0,\frac\pi2\bigr]$.
\end{thm}

As straightforward consequences of Theorem~\ref{Opeheim-Sin-Cos.tex-thm-1}, the following inequalities may be derived immediately.

\begin{thm}\label{Opeheim-Sin-Cos.tex-thm-2}
If $p\ge\frac2\pi$, then
\begin{equation}\label{Opeheim-Sin-Cos.tex-ineq}
\frac{(\pi/2)\sin x}{1+p\cos x}\le x\le\frac{(1+p)\sin x}{1+p\cos x},\quad 0\le x\le\frac\pi2;
\end{equation}
If $p\le\frac12$, the double inequality~\eqref{Opeheim-Sin-Cos.tex-ineq} reverses;
If $\frac12<p<\frac2\pi$, then
\begin{equation}\label{Opeheim-Sin-Cos.tex-ineq-min}
\frac{4p(1-p^2)\sin x}{1+p\cos x}\le x\le\frac{\max\{\pi/2, 1+p\}\sin x}{1+p\cos x}.
\end{equation}
The constants $\frac\pi2$ and $1+p$ in~\eqref{Opeheim-Sin-Cos.tex-ineq} and~\eqref{Opeheim-Sin-Cos.tex-ineq-min} are the best possible.
\end{thm}

\section{Concise proofs of theorems}

Now we are in position to concisely prove our theorems.

\begin{proof}[Proof of Theorem~\ref{Opeheim-Sin-Cos.tex-thm-1}]
Direct differentiation yields
\begin{align*}
f_p'(x)&=\frac{(x-\sin x\cos x)[p-(\sin x-x\cos x)/(x-\sin x\cos x)]}{(px\cos x +x)^2}\\
&\triangleq \frac{(x-\sin x\cos x)[p-h(x)]}{(px\cos x +x)^2},\\
h'(x)&=\frac{2\bigl[2 x^2+x\sin(2x)+2\cos(2x)-2\bigr]\sin x}{[2x-\sin(2x)]^2}\\
&\triangleq\frac{2g(x)\sin x}{[2x-\sin(2x)]^2},\\
g'(x)&=2 \cos (2 x) x+4 x-3 \sin (2 x),\\
g''(x)&=8(\tan x-x)\sin x\cos x\\
&>0
\end{align*}
on $\bigl(0,\frac\pi2\bigr)$. So the function $g'(x)$ is strictly increasing on $\bigl(0,\frac\pi2\bigr)$. Further, from $g'(0)=0$, it follows that $g'(x)>0$ and $g(x)$ is strictly increasing on $\bigl(0,\frac\pi2\bigr)$. Owing to $g(0)=0$, the functions $g(x)$ and $h'(x)$ is positive on $\bigl(0,\frac\pi2\bigr)$. As a result, the function $h(x)$ is strictly increasing on $\bigl(0,\frac\pi2\bigr)$. Due to $\lim_{x\to0^+}h(x)=\frac12$ and $h\bigl(\frac\pi2\bigr)=\frac2\pi$, it is concluded that
\begin{enumerate}
\item
when $p\ge\frac2\pi$, the derivative $f_p'(x)$ is positive on $\bigl(0,\frac\pi2\bigr)$, and so the function $f_p(x)$ is strictly increasing on $\bigl(0,\frac\pi2\bigr]$;
\item
when $p\le\frac12$, the derivative $f_p'(x)$ is negative on $\bigl(0,\frac\pi2\bigr)$, and so the function $f_p(x)$ is strictly decreasing on $\bigl(0,\frac\pi2\bigr]$;
\item
when $\frac12<p<\frac2\pi$, the derivative $f_p'(x)$ has a unique zero on $\bigl(0,\frac\pi2\bigr)$, and so the function $f_p(x)$ has a unique maximum on $\bigl(0,\frac\pi2\bigr]$.
\end{enumerate}
The proof of Theorem~\ref{Opeheim-Sin-Cos.tex-thm-1} is complete.
\end{proof}

\begin{proof}[Proof of Theorem~\ref{Opeheim-Sin-Cos.tex-thm-2}]
It is easy to see that
$$
\lim_{x\to0^+}f_p(x)=\frac1{1+p}\quad \text{and}\quad f_p\biggl(\frac\pi2\biggr)=\frac2\pi.
$$
By Theorem~\ref{Opeheim-Sin-Cos.tex-thm-1}, it follows that
\begin{enumerate}
  \item
  when $p\ge\frac2\pi$, we have
\begin{equation}\label{Opeheim-Sin-Cos.tex-pre}
  \frac1{1+p}<\frac{\sin x}{x(1+p\cos x)}\le \frac2\pi
\end{equation}
on $\bigl(0,\frac\pi2\bigr]$, which may be rewritten as the inequality~\eqref{Opeheim-Sin-Cos.tex-ineq};
  \item
  when $p\le\frac12$, the inequality~\eqref{Opeheim-Sin-Cos.tex-pre} reverses;
  \item
  when $\frac12<p<\frac2\pi$, we have
$$
\frac{\sin x}{x(1+p\cos x)}>\min\biggl\{\frac1{1+p}, \frac2\pi\biggr\}
$$
on $\bigl(0,\frac\pi2\bigr)$, which may be rearranged as the right-hand side inequality in~\eqref{Opeheim-Sin-Cos.tex-ineq-min}.
\end{enumerate}
\par
The left-hand side inequality in~\eqref{Opeheim-Sin-Cos.tex-ineq-min} can be deduced by the same argument as in~\cite[p.~60]{Oppeheim-Zhu-mia}.
The proof of Theorem~\ref{Opeheim-Sin-Cos.tex-thm-2} is complete.
\end{proof}

\section{Remarks}

After proving our theorems, we giver several remarks on them.

\begin{rem}
For $p\le\frac12$, the reversed version of the inequality~\eqref{Opeheim-Sin-Cos.tex-ineq} may be rewritten as
\begin{equation}\label{Opeheim-Sin-Cos.tex-rew-1}
\frac{2(1+p\cos x)}\pi<\frac{\sin x}{x}\le \frac{(1+p\cos x)}{1+p},\quad 0<x\le\frac\pi2.
\end{equation}
Integrating on both sides of~\eqref{Opeheim-Sin-Cos.tex-rew-1} gives
\begin{equation*}
1+\frac2{\pi}p<\int_0^{\pi/2}\frac{\sin x}{x}\td x< \frac{\pi+2p}{1+p},\quad p\le\frac12.
\end{equation*}
Hence, taking $p=\frac12$ in the above inequality leads to
\begin{equation}\label{int-sin-x-ineq-1}
1.31\dotsm=1+\frac1{\pi}<\int_0^{\pi/2}\frac{\sin x}{x}\td x<\frac{2(\pi+1)}{3}=2.76\dotsm.
\end{equation}
Similarly, if integrating and letting $p=\frac2\pi$ in~\eqref{Opeheim-Sin-Cos.tex-ineq}, then
\begin{equation}\label{int-sin-x-ineq-2}
1.34\dotsm=\frac{4+\pi^2}{2(2+\pi)}<\int_0^{\pi/2}\frac{\sin x}{x}\td x <1+\biggl(\frac2{\pi}\biggr)^2=1.40\dotsm.
\end{equation}
\end{rem}

\begin{rem}
For $\frac12<p<\frac2\pi$, the inequality~\eqref{Opeheim-Sin-Cos.tex-ineq-min} may be rearranged as
\begin{equation}\label{Opeheim-Sin-Cos.tex-ineq-min-rew}
\min\biggl\{\frac2\pi, \frac1{1+p}\biggr\}(1+p\cos x)\le \frac{\sin x}{x}\le \frac{1+p\cos x}{4p(1-p^2)},\quad 0<x\le\frac\pi2.
\end{equation}
As done above, integrating gives
$$
\min\biggl\{\frac2\pi, \frac1{1+p}\biggr\}\biggl(p+\frac\pi2\biggr)<\int_0^{\pi/2}\frac{\sin x}{x}\td x <\frac{2p+\pi}{8p\bigl(1-p^2\bigr)},\quad \frac12<p<\frac2\pi.
$$
Maximizing the lower bound and minimizing the upper bound in the above double inequality reduce to
\begin{equation}\label{int-sin-x-ineq-3}
1.36\dotsm=2\biggl(1-\frac1\pi\biggr)<\int_0^{\pi/2}\frac{\sin x}{x}\td x <\frac{2p_0+\pi}{8p_0\bigl(1-p_0^2\bigr)}=1.37\dotsm,
\end{equation}
where
\begin{align*}
p_0&=\frac{\pi}{4}\biggl\{\cos \biggl[\frac{1}{3} \arctan\biggl(\frac{4\sqrt{\pi^2-4}\,} {\pi^2-8}\biggr)\biggr] +\sqrt{3}\,\sin\biggl[\frac{1}{3} \arctan\biggl(\frac{4\sqrt{\pi^2-4}\,} {\pi^2-8}\biggr)\biggr]-1\biggr\}\\
&=0.52\dotsm.
\end{align*}
\par
Comparing the inequalities~\eqref{int-sin-x-ineq-1}, \eqref{int-sin-x-ineq-2} and~\eqref{int-sin-x-ineq-3} shows that the inequality~\eqref{Opeheim-Sin-Cos.tex-ineq-min} or~\eqref{Opeheim-Sin-Cos.tex-ineq-min-rew} is more accurate in whole.
\par
The inequality~\eqref{int-sin-x-ineq-3} improves the inequalities
\begin{equation}
1.33\dotsm=\frac43<\int_0^{\pi/2}\frac{\sin x}{x}\td x<\frac{\pi+1}3=1.38\dotsm
\end{equation}
and
\begin{equation}
\int_0^{\pi/2}\frac{\sin x}{x}\td x>\frac{\pi+5}6=1.35\dotsm
\end{equation}
obtained in~\cite[p.~521, (32)]{construct} and~\cite{qi-jordan}.
\end{rem}

\begin{rem}
In~\cite[p.~247, 3.4.31]{mit}, it was listed that the inequality
\begin{equation}\label{Shafer-ineq-arcsin}
\arcsin x>\frac{6\bigl(\sqrt{1+x}\,-\sqrt{1-x}\,\bigr)}{4+\sqrt{1+x}\,+\sqrt{1-x}\,} >\frac{3x}{2+\sqrt{1-x^2}\,}
\end{equation}
holds for $0<x<1$. It was also pointed out in~\cite[p.~247, 3.4.31]{mit} that these inequalities are due to R. E. Shafer, but no a related reference is cited. By now we do not know the very original source of inequalities in~\eqref{Shafer-ineq-arcsin}.
\par
In the first part of the short paper~\cite{Fink-Beograd-Univ-95}, the inequality between the very ends of~\eqref{Shafer-ineq-arcsin} was recovered and an upper bound for the arc sine function was also established as follows:
\begin{equation}\label{Shafer-Fink-ineq-arcsin}
\frac{3x}{2+\sqrt{1-x^2}\,}\le\arcsin x\le\frac{\pi x}{2+\sqrt{1-x^2}\,},\quad 0\le x\le 1.
\end{equation}
Therefore, we call~\eqref{Shafer-Fink-ineq-arcsin} the Shafer-Fink's double inequality for the arc sine function.
\par
In~\cite{Malesevic-Beograd-Univ-97}, the right-hand side inequality in~\eqref{Shafer-Fink-ineq-arcsin} was improved to
\begin{equation}\label{Malesevic-Beograd-Univ-97-ineq}
\arcsin x\le\frac{\pi x/(\pi-2)}{2/(\pi-2)+\sqrt{1-x^2}\,},\quad 0\le x\le1.
\end{equation}
\par
As done in~\cite{Oppeheim-Zhu-mia}, by taking $t=\sin x$ in Theorem~\ref{Opeheim-Sin-Cos.tex-thm-2}, the inequalities in~\eqref{Shafer-Fink-ineq-arcsin}, \eqref{Malesevic-Beograd-Univ-97-ineq} and the following Shafer-Fink type inequalities may be derived readily:
\begin{gather}\label{Zhu-mia-1}
\frac{\pi(4-\pi)x}{2/(\pi-2)+\sqrt{1-x^2}\,} \le \arcsin x,\quad 0\le x\le1;\\
\frac{\pi x/2}{1+\sqrt{1-x^2}\,} \le\arcsin x,\quad 0\le x\le1.\label{Zhu-mia-3}
\end{gather}
\par
All corresponding bounds in~\eqref{Shafer-Fink-ineq-arcsin}, \eqref{Malesevic-Beograd-Univ-97-ineq}, \eqref{Zhu-mia-1} and~\eqref{Zhu-mia-3} are not included each other.
\par
The facts stated above strongly shows us that Oppenheim type inequalities and Shafer-Fink type inequalities can be converted to each other.
\end{rem}

\begin{rem}
Recently, some new Shafer-Fink type inequalities and generalizations of Oppenheim's inequality are procured in~\cite{Baricz-Zhu-JIA, Shafer-Fink-ArcSin.tex, Zhu-New-Arc-Hyper, Zhu-Pan-Shafer-Fink}.
\end{rem}

\begin{rem}
In~\cite{AVV-Pacific-93, Pinelis-JIPAM-02-01}, the following L'H\^ospital rule for monotonicity was established:
Let $f(x)$ and $g(x)$ be continuous functions on $[a,b]$ and differentiable on $(a,b)$ such that $g'(x)\ne0$ on $(a,b)$. If $\frac{f'(x)}{g'(x)}$ is increasing \textup{(}or decreasing respectively\textup{)} on $(a,b)$, then the functions $\frac{f(x)-f(b)}{g(x)-g(b)}$ and $\frac{f(x)-f(a)}{g(x)-g(a)}$ are also increasing \textup{(}or decreasing respectively\textup{)} on $(a,b)$.
This conclusion has been employed in a lot of literature such as~\cite{jordan-strengthened, jordan-generalized.tex} and closely-related references therein. This conclusion can also be utilized to prove the increasing monotonicity of the function $h(x)$ in the proof of Theorem~\ref{Opeheim-Sin-Cos.tex-thm-1} as follows.
\par
Let $h_1(x)=\sin x-x\cos x$ and $h_2(x)=x-\sin x\cos x$ on $\bigl[0,\frac\pi2\bigr]$. Then
$$
h_1'(x)=x\sin x, \quad h_2'(x)=2\sin^2x,
$$
and so
$$
\frac{h_1'(x)}{h_2'(x)}=\frac{x}{2\sin x}
$$
is strictly increasing on $\bigl(0,\frac\pi2\bigr)$. Consequently, the function
$$
h(x)=\frac{h_1(x)}{h_2(x)}=\frac{h_1(x)-h_1(0)}{h_2(x)-h_2(0)}
$$
is strictly increasing on $\bigl(0,\frac\pi2\bigr)$.
\end{rem}

\begin{rem}
It is worthwhile to noting that the surname name ``Oppenheim'' was mistaken for ``Oppeheim'' in~\cite{Oppeheim-Zhu-mia}.
\end{rem}

\end{document}